\documentclass[12pt]{article}
\usepackage{amsmath,amsthm,amsfonts,amssymb}
\newtheorem{theorem}{Theorem}
\newtheorem{lemma}{Lemma}
\newtheorem{cor}{Corollary}

\newtheorem{prop}{Proposition}
\newcommand{\lla}{\langle\!\langle}
\newcommand{\rra}{\rangle\!\rangle}
\renewcommand{\lll}{|\!|\!|}
\newcommand{\LASEP}{L}
\newcommand{\ASEP}{\lll_{-1,\lambda}}

\newcommand{\Z}{\mathbb{Z}}

\newcommand{\la}{\langle}
\newcommand{\ra}{\rangle}

\begin{document}

\title{$t^{1/3}$ Superdiffusivity of Finite-Range Asymmetric Exclusion Processes on $\mathbb Z$ }
\author{Jeremy Quastel, Benedek Valk\'o
\\University of Toronto}
\maketitle

\footnotetext[1]{Supported by the Natural Sciences and Engineering Research Council of Canada.

{\em Key words and phrases.}  asymmetric exclusion process, superdiffusivity. 

{\em MSC}: 60K35, 82C22.}




\abstract{We consider finite-range asymmetric exclusion processes on $\mathbb Z$
with non-zero drift.  The diffusivity $D(t)$ is expected to
be of ${\mathcal O}(t^{1/3})$.  We prove that  $D(t)\ge Ct^{1/3}$ in the weak (Tauberian)
sense that $\int_0^\infty e^{-\lambda t }tD(t)dt \ge 
C\lambda^{-7/3}$ as $\lambda\to 0$.   The proof employs the resolvent method to make a direct comparison with
the totally asymmetric simple exclusion process, for which the result is a consequence
of the scaling limit for the two-point function recently obtained by Ferrari and
Spohn.  
In the nearest neighbor case, we show further that $tD(t)$ is monotone, and hence we can conclude that $D(t)\ge Ct^{1/3}(\log t)^{-7/3}$ in the usual sense.}

\section{Introduction}
\label{intro}
\setcounter{equation}{0}

A finite-range exclusion process on the integer lattice $\mathbb Z$ is a system of 
continuous time, rate one random walks with finite-range jump law $p(\cdot)$, i.e. $p(z)\ge 0$, and $p(z)=0$
for $z>R$ for some $R<\infty$, $\sum_z p(z)=1$, interacting via {\it exclusion}:  Attempted jumps to occupied
sites are suppressed.  We will always assume in this article that $p(\cdot)$ has a 
non-zero drift, \begin{equation}\sum_z zp(z) = b\neq 0.\label{drift}\end{equation}  In particular, $p(\cdot)$ is asymmetric and we will refer to the process as the asymmetric exclusion process (AEP).  The state space of the process is $\{0,1\}^\mathbb Z$ and it is traditional to call
configurations $\eta$, with $\eta_x\in\{0,1\}$ indicating the absence, or presence, of a 
particle at $x\in \mathbb Z$. The infinitesimal generator of the process is given by
\begin{equation}
\LASEP f(\eta)=
\sum_{x,z\in \mathbb Z} p(z) \eta_x (1-\eta_{x+z})(f(\eta^{x,x+z}) -f(\eta) ) 
\end{equation}
where $\eta^{x,y}$ denotes the configuration obtained from $\eta$ by interchanging the
occupation variables at $x$ and $y$.  The Bernoulli product measures $\pi_\rho$, $\rho\in [0,1]$, with $\pi_\rho(\eta_x=1) = \rho$ form a one-parameter family of invariant measures for the process. 
In fact, there exist other invariant measures \cite{BM}, but they
will not be relevant for our discussion. The process starting from $\pi_0$ and $\pi_1$ are trivial and so we
consider the stationary process obtained by starting with $\pi_\rho$ for some $\rho\in (0,1)$.  

Let
\begin{equation}
\hat\eta_x = \frac{ \eta_x-\rho}{ \sqrt{\rho(1-\rho)}}, \qquad
\hat\eta_{A} = \prod_{x\in A} \hat\eta_x
\end{equation}
for any finite nonempty set $A\subset\Z$.  The collection $\{\hat \eta_A\}$ where $A$ ranges over finite subset of $\Z$ is an orthonormal basis of $L^2(\pi_\rho)$ with its natural inner product
\begin{equation}
\langle f,g\rangle = \int_{\{0,1\}^\mathbb Z} fg d\pi_\rho.
\end{equation}
  Then
$L^2(\pi_\rho)$ can naturally be thought of as the direct sum of subspaces $H_1, H_2,\ldots$
where $H_n$ is the linear span of $\{\hat \eta_A\}, ~ |A|=n$.   It is natural to think of 
$H_1$ as being {\it linear } functions, $H_2$ as {\it quadratic functions}, etc.

From a physical point of view, the most basic quantity is the two-point function,
\begin{equation}
S(x,t) = E[ (\eta_x(t)-\rho)( \eta_0(0)-\rho)].
\end{equation}
The expectation is with respect to the stationary process obtained by starting from one of the invariant measures
$\pi_\rho$.
It is easy to show (see \cite{PS}) that $S(x,t)$ satisfies the sum rules
\begin{equation} \label{MV}
\sum_x S(x,t) = \rho(1-\rho) = \chi,\qquad \frac{1}{\chi}\sum_x x S(x,t) = (1-2\rho)bt.
\end{equation}
Note that one should not expect to be able to actually compute $S(x,t)$ but one does hope to find
its large scale structure.  The next most basic quantity, the diffusivity
$D(t)$,  is already unknown.  It is defined as
\begin{equation}\label{dee}
D(t) = ({\chi t})^{-1} \sum_{x\in \mathbb Z} (x - 
(1-2\rho) b t )^2 S(x,t).
\end{equation}

Using coupling (see \cite{L}), the diffusivity can be rewritten in
terms of the variance of a second class particle.
Suppose  one starts with two configurations $\eta'$ and $\eta$ which are ordered in the
sense that $\eta'_x\ge \eta_x$ for each $x\in \mathbb Z$.  One can couple the two exclusions
by having them jump together whenever possible and one observes that at later times
the ordering is preserved.  If we write $\eta' = \eta + \eta''$ then the "particles"
of $\eta''$ move according to the {\it second class } particle dynamics.  Among themselves
they move with the standard exclusion rule, the other (first class) particles move without noticing them,
and if a first class particle attempts to jump to a site occupied by a second class particle,
the two exchange positions.  Note that
\begin{eqnarray}\label{gra}
\chi^{-1} S(x,t)  & =  & P( \eta_x(t)=1 ~|~\eta_0(0)=1)- P(\eta_x(t)=1~|~\eta_0(0)=0)
 \nonumber \\& = &  {\mathbb P}( \eta''_x(t) =1~|~ \eta''(0) =\delta_0) 
 \nonumber\\ & = & P(X(t)=x~|~ X(0)=0) 
\end{eqnarray}
are the transition probabilities of a single second class particle $X(t)$ starting at the origin.  Here $\delta_0$ is the configuration with only one particle at $0$ and $\mathbb P$ is the coupled measure. The diffusivity is then given by\begin{equation}
D(t)= t^{-1} Var (X(t)).
\end{equation}

 We can alternately write the dynamics as a stochastic differential equation
 \begin{equation}
 {\mathbf d} \hat\eta_x = (\nabla w_x + \Delta\hat\eta_x ) dt + d{\mathcal M}_x
 \end{equation}
 where $\mathbf d$ is a microscopic convective derivative,
 \begin{equation}
 \mathbf d\hat\eta_x = d\hat \eta_x +  \frac {(1-2\rho)}{2}\sum_zp(z) (\hat\eta_{x+z} -\hat\eta_{x-z})
 \end{equation}
  $\nabla$ and $\Delta$ are microscopic analogues of first and second
 spatial derivatives,
 \begin{eqnarray}
 \nabla w_x & = &\chi^{1/2}\sum_zp(z) (\hat\eta_{x,x+z}- \hat\eta_{x-z,x}),\\
\Delta\hat\eta_x & = &\frac{1}{2}\sum_z p(z) (\hat\eta_{x+z} +\hat\eta_{x-z} -2\hat\eta_x),
 \end{eqnarray}
 and ${\mathcal M}_x(t)$ are martingales with
\begin{equation}
E[ ( \int_0^t \sum_x \phi_x d{\mathcal M}_x )^2 ] = \int_0^t \sum_{x,z} p(z) 
(\phi_{x+z}-\phi_x)^2 ds.
\end{equation}
The current
$w_x= \tau_x w$, where the specific quadratic function $w$ is given by
\begin{equation}\label{w}
w =\sum_z z p(z) \hat \eta_{\{0,z\}}.
\end{equation}

In this sense, AEP is a natural discretisation of the stochastic Burgers equation,
\begin{equation}\label{SBE}
\partial_t u = \partial_x u^2 + \partial_x^2 u + \partial_x \dot W
\end{equation}
for a function $u(t,x)$ of $x\in \mathbb R$ and $t>0$ where $\dot W$ is a space-time
white noise.  White noise is supposed to be an invariant measure.  Letting $\partial_x U = u$ one obtains the Kardar-Parisi-Zhang equation
for surface growth,
\begin{equation}
\partial_t U = (\partial_xU)^2 + \partial_x^2 U + \dot W.
\label{KPZ}
\end{equation}
We are interested in the large scale behaviour and 
the only rescalings of $u$ which preserve the initial white noise are
\begin{equation}
u_\epsilon(t,x) = \epsilon^{-1/2} u (\epsilon^{-z} t,\epsilon^{-1} x).
\end{equation}
The stochastic Burgers equation (\ref{SBE}) transforms to,
\begin{equation}
\partial_t u_\epsilon = \epsilon^{ \frac{3}{2} - z} \partial_x u_\epsilon^2 + \epsilon^{2-z}\partial_x^2 u_\epsilon + \epsilon^{1-\frac{z}{2}} \partial_x\dot W ,
\end{equation}
which suggests that the dynamical exponent $z={3}/{2}$ and that the diffusion and 
random forcing terms become irrelevant in the limit. 

 The exponent $z={3}/{2}$ was 
first predicted for (\ref{SBE}) by Forster, Nelson and Stephen \cite{FNS}, then for 
AEP by van Beijeren, Kutner and Spohn \cite{BKS} and then for
(\ref{KPZ}) by Kardar, Parisi and Zhang \cite{KPZ}.  Note that at a rigorous level we are
very far from  understanding this for either (\ref{SBE}) or (\ref{KPZ}).
At the present time the mathematical
 problem there is just to make sense of the equation (see \cite{BG}).
So it makes sense to consider exclusion processes, which are clearly 
well defined, yet are supposed to have the same large scale behaviour.

The scaling prediction for $u$ suggests that on  large scales
\begin{equation}
S(x,t) \simeq  t^{-2/3} \Phi( t^{-2/3}(x- (1-2\rho) bt))
\end{equation}
for some  scaling function $\Phi$, and in particular one conjectures that,
\begin{equation}\label{conj}
D(t) \simeq C t^{1/3}.
\end{equation}

Note that the case of asymmetric exclusion with mean-zero jump law is different and there one
has as usual that $D(t)\to D$ as $t\to\infty$ (see \cite{V}).

The diffusivity can be related to the time integral
of current-current correlation functions by the Green-Kubo formula,
\begin{equation}
D(t) = \sum_z z^2 p(z)  +2\chi t^{-1} \int_0^t\int_0^s \lla w, e^{uL} w\rra duds. \label{GK}
\end{equation}
 It uses a 
  special inner product defined for local
functions  by
\begin{equation}
\lla \phi, \psi\rra =\langle \phi, \sum_x\tau_x \psi\rangle.
\end{equation}
(\ref{GK}) is proved in \cite{LOY} (in the special case $p(1)=1$,
but the proof for general AEP is the same.) 
A useful variant is obtained by taking the Laplace transform,
\begin{equation}
\int_0^\infty e^{-\lambda t} tD(t) dt = \lambda^{-2} \left(
\sum_z z^2 p(z)+ 2\chi \lll w\ASEP^2 \right)
\end{equation}
where the $H_{-1}$ norm corresponding to $\LASEP$ is defined for local
functions by
\begin{equation}
\lll\phi\ASEP = \lla \phi, (\lambda - \LASEP)^{-1} \phi \rra^{1/2} .
\end{equation}

We say that $D(t) \simeq t^{\rho}$, $\rho>0$ {\it in the weak (Tauberian) sense}
if $\int_0^\infty e^{-\lambda t} tD(t) dt\simeq \lambda^{-(2 + \rho)}$.  Hence
the weak (Tauberian) version of the conjecture (\ref{conj}) is
\begin{equation}
\lll w\ASEP^2 \simeq \lambda^{-1/3}.
\end{equation}
One of the key advantages of this {\it resolvent approach} is that there is a variational formula (see \cite{LQSY}),
\begin{equation}\label{varform}
\lll w\ASEP^2 = \sup_f \left\{ 2\lla w,f\rra- \lla f,(\lambda-S) f\rra-\lla Af, (\lambda-
S)^{-1} Af\rra \right\}
\end{equation}
where $S= \frac{1}{2} (L+L^*)$ and $A=\frac{1}{2} (L-L^*)$ are the symmetric and antisymmetric
parts of the generator $L$.  $S$ is nothing but the generator of the symmetric exclusion
process with $\bar p(z) =\frac12(p(z)+p(-z))$.  It has the special property that it maps the subspaces $H_n$ into themselves,
and on each is nothing but the generator of a symmetric random walk.  Hence one can hope
to obtain non-trivial information from (\ref{varform}) by choosing carefully test functions
$f$.  This idea was used  in \cite{LQSY} to obtain $D(t) \ge C t^{1/4}$ in $d=1$ and $D(t)
\ge C(\log t)^{1/2}$ in $d=2$, which was improved to $D(t)\simeq C(\log t)^{2/3}$ in \cite{Y}.  All of these are in the weak (Tauberian) sense.  

The special case of jump law $p(1)=1$, $p(z)=0$, $z\neq1$ is called the totally asymmetric
simple exclusion process (TASEP).  Simple refers to the nearest-neighbour jumps of the 
underlying random walk.  It is very remarkable that after about 20 years of intense study,
TASEP has succumbed to a combination of sophisticated techniques from analysis, combinatorics
and random matrix theory (see \cite{FS} and references therein).

We now state the main result of Ferrari and Spohn \cite{FS}.  Define the height function
\begin{equation}\label{h}
h_t(x) =   2N_t -M_t(x)
\end{equation}
$t\ge 0$, where $N_t$ counts the number of jumps from site $0$ to site $1$ up to time $t$ and
\begin{equation}\label{def_M}
M_t(x)=\left\{
\begin{array}{lc}
  \phantom{-}\sum_{i=1}^x (2\eta_i(t)-1)&   \textup{if } x>0,   \\[7pt]
  \phantom{MMM}0& \textup{if } x=0\\[7pt]
  -\sum_{i=x+1}^0 (2\eta_i(t)-1)&   \textup{if } x<0. 
\end{array}
 \right.
\end{equation} Note
that $
E[ h_t(x) ] = 2\chi t + (1-2\rho) x$.   Let 
\begin{equation}\label{veee}
v(x,t) = Var( h_t (x) ).
\end{equation}
Since $h_t(x+1) -h_t(x) =1-2\eta_{x+1}(t)$,  it is not hard to check  that 
\begin{equation}\label{stov}
8S(x,t) = v(x+1,t)-2v(x,t) + v(x-1,t).
\end{equation} See \cite{PS} for a detailed proof.  We have
\begin{equation}\label{d-var}
D(t) = (4\chi t)^{-1} \sum_{x\in \mathbb Z} Var( h_t (x) )-4\chi |x-(1-2\rho)t| .
\end{equation}
(see Section \ref{section_ps}.)
Now consider a normalised version of $h_t$:
\begin{equation}
\hat h_t (x) = \chi^{-2/3} t^{-1/3} ( h_t(x) - E[ h_t(x)]).
\end{equation}
and for each fixed $t>0$ and $\omega\in\mathbb R$ let $F_{\omega,t}$ be
the cumulative distribution function of $-\hat h_t(\lfloor (1-2\rho) t + 2\omega\chi^{1/3} t^{2/3}\rfloor)$;
\begin{equation}
F_{\omega,t}(s) = P(- \hat h_t(\lfloor (1-2\rho) t + 2\omega\chi^{1/3} t^{2/3}\rfloor)\le s)
\end{equation}

The main result of Ferrari and Spohn concerning  TASEP  is
that 
 $dF_{\omega,t}$ converge weakly as probability measures, as $t$ tends to infinity,
to $dF_\omega$ where
\begin{equation}
F_\omega(s) = \frac{\partial}{\partial s} \left( F_{GUE} (s+\omega^2)
g(s+\omega^2,\omega) \right)
\end{equation}
 where $F_{GUE}$ is the Tracy-Widom distribution
and $g$ is a scaling function defined through the Airy kernel (see \cite{FS} for details).

Note that the convergence stated in \cite{FS} is that
for any $c_1<c_2$, 
\begin{equation} \label{FSmain}
\lim_{t\to\infty} \int_{c_1}^{c_2} F_\omega(s,t) ds  = \int_{c_1}^{c_2}F_\omega(s)ds.
\end{equation}
In fact, this is the same as weak convergence.
For by monotonicity, if  $\epsilon>0$,
\begin{equation}
 \epsilon^{-1} \int_{s-\epsilon}^s F_{\omega,t}(u)du
\le F_{\omega,t}(s)\le \epsilon^{-1} \int_s^{s+\epsilon} F_{\omega,t}(u)du.
\end{equation}
Taking the limit in $t$ and using (\ref{FSmain}) we see that $\lim_{t\to\infty} F_{\omega,t}(s) = F_{\omega}(s)$ at any continuity point $s$ of the
limit function (in this case all $s\in \mathbb R$), and this is equivalent to weak convergence.
 
The proof of Ferrari and Spohn is through a direct mapping between TASEP and a particular last passage percolation
problem.  Such a mapping is not available except for the case of TASEP.
So although one expects analogous results for general AEP in one
dimension, different techniques will be required.  Our main motivation
here is to confirm, at least in part, the predicted universality
(see Section 6 of \cite{PS} for a nice description) by showing
that these results for TASEP imply some bounds for general AEP.

From (\ref{d-var}) and (\ref{FSmain}) one expects 
\begin{equation}
D^{TASEP}(t)\simeq c^{TASEP} \chi^{2/3} t^{1/3}
\end{equation}
where 
\begin{equation}
c^{TASEP} = \int d\omega\int s^2 dF_\omega(s)= 2\int d\omega \int ds  F_{GUE} (s+\omega^2)
g(s+\omega^2,\omega) .
\end{equation} Here, and throughout this article, we will use the superscript $\hbox{}^{TASEP}$ to denote the values taken by TASEP of quantities defined for general AEP.
Unfortunately, the necessary estimates for the upper bound appear to be missing at this time.   
However from the weak convergence we have immediately that  
\begin{cor}\label{cor1}
\begin{equation}
\liminf_{t\to\infty } t^{-1/3} D^{TASEP}(t) \ge c^{TASEP} \chi^{2/3}.\label{lbtasep}
\end{equation}
\end{cor}

{\it Remark.} Another way to see the strict positivity 
of the left hand side without computing 
$c^{TASEP}$
is that by Schwartz's inequality and (\ref{dee}),
\begin{equation}
D(t) \ge t^{-1}\left({\chi }^{-1} \sum_{x\in \mathbb Z} |x - 
(1-2\rho)  t| S(x,t)\right)^2 .
\end{equation}
We have
\begin{equation}\label{abs_mom}
\sum_{x\in \mathbb Z} |x - 
\lfloor(1-2\rho)  t\rfloor | S(x,t)= 2Var(h_t(\lfloor (1-2\rho) t \rfloor))
\end{equation} (see Section \ref{section_ps})
and from the weak convergence we have,
\begin{equation}
\liminf_{t\to\infty} t^{-2/3}Var(h_t(\lfloor (1-2\rho) t \rfloor)) 
\ge \chi^{4/3}\int s^2 dF_\omega (s).
\end{equation}
Since by (\ref{MV})
\begin{equation*}
\left| \sum_{x\in \mathbb Z} \left(|x - 
\lfloor(1-2\rho)  t\rfloor |- |x - 
(1-2\rho)  t |\right)S(x,t)\right|\le\chi
\end{equation*}
the positive lower bound on $\liminf_{t\to\infty } t^{-1/3} D^{TASEP}(t) $ follows.\bigskip

The main result of the present article  is a comparison between the diffusivity of 
AEP and that of TASEP:

\begin{theorem}\label{thm_main}
Let $D(t)$ be the diffusivity of a 
 finite range exclusion process in $d=1$
with non-zero drift.  Let $D^{TASEP}(t)$ be the 
diffusivity of the totally asymmetric simple exclusion process.  There exists $C<\infty$ such that 
\begin{eqnarray}\label{equiv}
C^{-1}\int_0^\infty e^{-\lambda t} t D^{TASEP}(t) dt
&\le& \int_0^\infty e^{-\lambda t} t D(t) dt \\ \nonumber &
\le & C\int_0^\infty e^{-\lambda t} t D^{TASEP}(t) dt
\end{eqnarray}\end{theorem}

Combined with   (\ref{lbtasep}) this gives
\begin{theorem} \label{thm_tauber}For any finite range exclusion process in $d=1$
with non-zero drift,
$D(t)\ge Ct^{1/3}$ in the weak (Tauberian) sense:
There exists $C>0$ such that
\begin{equation}\label{l_bound}
\int_0^\infty e^{-\lambda t} t D(t) dt \ge C\lambda^{-7/3}.
\end{equation}
\end{theorem}

 We now make some comments on obtaining strict versions of the
estimates, as opposed to weak (Tauberian) versions. 

 In \cite{LY}
it is shown that
\begin{equation}
t^{-1}\sum_x E[\int_0^t w (s) ds \int_0^t \tau_x w (s) ds] \le \lll w\lll_{-1, t^{-1}}^2
\end{equation}
and hence an upper weak (Tauberian) bound implies a strict upper bound
in time on the diffusivity.  There is no analogous fact for lower bounds.  However, it
is easy to show the following: 
\begin{prop}\label{tauber}
Suppose that  $v(t)\ge 0$ is a nondecreasing function and $\beta>0$.

\noindent 1.  Suppose there exist $c_1<\infty$
and $\lambda_0>0$ such that for $0<\lambda<\lambda_0$,
\begin{equation}
 \int_0^\infty e^{-\lambda t} v(t) dt \le c_1\lambda^{-(1+\beta)} 
\end{equation}
 then there exist $c_2<\infty$ and $t_0$ such that for all
$t>t_0$,
\begin{equation}
 v(t) \le c_2 t^\beta.
\end{equation}
2.  Suppose $v(t) \le c_2t^\alpha$ for some $\alpha\ge\beta$ and $t>t_0$
and for some $c_3>0$, for $0<\lambda<\lambda_0$,
\begin{equation}
\int_0^\infty e^{-\lambda t} v(t) dt \ge c_3\lambda^{-(1+\beta)}.
\end{equation} 
Then there exists $c_4>0$ and $t_1<\infty$ such that for $t>t_1$
\begin{equation}
v(t) \ge  
\begin{cases} 
c_4 t^\beta & {\rm if}~  \alpha=\beta;\\
c_4t^{\beta} (\log t)^{-(1+\beta)} & {\rm if}~  \alpha>\beta.
\end{cases}
\end{equation}
\end{prop}
\begin{proof} 1. Since $v$ is monotone nondecreasing we have
for $t>\lambda_0^{-1}$,
\begin{equation}
e^{-1}v(t)= \int_1^\infty e^{-s} v(t) ds \le \int_0^\infty e^{-s} v(t s) ds\le  c_1 t^{\beta}.
\end{equation} 
 2. Because $v(t)$ is non-decreasing, $\int_0^t e^{-\lambda s}v(s) ds \le tv(t)$ and if $v(t)\le c_2t^\alpha$ we have $\int_t^\infty e^{-\lambda s} v(s)ds \le c_2' \lambda^{-1} e^{-\lambda t} t^\alpha$ for $t>t_1$.
Hence
\begin{equation}
c_3\lambda^{-(1+\beta)}\le tv(t) + c_2'  \lambda^{-1} e^{-\lambda t} t^\alpha.
\end{equation}
Choosing $\lambda = t^{-1}( 1+ (\alpha-\beta) (\log t+ c\log\log t))$ gives the result.  
\end{proof}

Note that the bound 
\begin{equation} \int_0^\infty e^{-\lambda t} tD(t) dt \le C\lambda^{-5/2}
\end{equation}
can be derived easily from the variational formula (\ref{varform}) (see the proof of Proposition 
\ref{prop2} for a similar computation).

Certainly one expects $tD(t)$ to be nondecreasing in general.  
We will show in Lemma \ref{lemma2} that\begin{equation}
\partial_t (tD(t)) = \sum_z z^2p(z) -2  \rho \sum_{z>0}  z (p(z)-p(-z)) E[
\widetilde X(t)|\eta_{z}(0)=1]
\end{equation}
where
\begin{equation}
\widetilde X(t) = X(t)-(1-2\rho) b t,
\end{equation}
What one expects is that  $b\, E[\widetilde X(t)~|~ \eta_{z}(0) =1]\le 0$.  
If $p(z)\ge p(-z)$ for all $z>0$, (or for all $z<0$) this would
imply that $tD(t)$ is increasing.  We have only been able to prove
this in the special case of the simple (nearest neighbor) exclusion
(see Proposition \ref{prop:mon}).  
Hence for this class of AEP we can make the following statement:

\begin{theorem} Let $D(t)$ be the diffusivity of a nearest neighbor ($p(z)=0$, $|z|\neq 1$) asymmetric exclusion. 

\noindent 1. There exists $c_0>0$ such that
\begin{equation}
D(t) \ge c_0 t^{1/3} (\log t)^{-7/3}.
\end{equation}
2. Suppose that there exists $c_1<\infty$ such that 
\begin{equation}
D^{TASEP}(t) \le c_1 t^{1/3}.
\end{equation} Then there exists $c_2<\infty$ such that
\begin{equation}
c_2^{-1} t^{1/3} \le D(t) \le
c_2 t^{1/3}.
\end{equation}
\end{theorem}
{\it Remarks.} 1.  Note that in Theorems \ref{thm_main} and \ref{thm_tauber} we have not made any assumptions about the
irreducibility of $p(\cdot)$. Let \begin{equation}
\kappa=\textup{gcd}(y\in \Z: p(y)>0).
\end{equation}
If $\kappa>1$ then our AEP is the same as $\kappa$ independent copies of 
the AEP with jump law $\tilde p(y)=p(\kappa y)$ on the sublattices $\kappa 
\Z +i$ ($i=0,1,\dots,\kappa-1$). Using this simple observation it is easy 
to extend all our proofs from $\kappa=1$ to $\kappa>1$, so we can assume
without loss of generality in the proofs that $p(\cdot)$ is irreducible.

{2.} Analogous methods to the ones presented here could in principle
be applied to other functionals of AEP.  For example, 
the variance of the  occupation time of the origin, 
\begin{equation}
\int_0^t \eta_s(0) ds,
\end{equation}
 is also expected to be $\mathcal O(t^{4/3})$.  In \cite{B} a lower
 bound of the form $Ct^{5/4}$ is obtained.  This variance is again given
 by the $H_{-1}$ norm of a certain function and direct comparisons between
  its value for TASEP and general AEP can be obtained in a straightforward way.  Hence asymptotic order of growth 
  bounds for this variance under TASEP
  would imply the same for AEP.  Unfortunately, at the present time no
  such bounds are available, though it is plausible they could
  be derived from the machinery that has been developed for TASEP.

\section{Comparison of $H_{-1}$ norms}
\label{}
\setcounter{equation}{0}


The first proposition adapts  results of Sethuraman \cite{S} to
the present context.
 
\begin{prop}  \label{prop1} There exist $\alpha, \beta\in (0,\infty)$ depending
only on $p(\cdot)$ such that
\begin{equation}
\alpha^{-1}\lll \phi\lll_{-1, \beta^{-1}\lambda}^{TASEP}  \le
\lll \phi \lll_{-1,\lambda} \le \alpha
\lll \phi\lll_{-1,\beta\lambda}^{TASEP} 
\end{equation}
\end{prop}

\begin{proof}
We can also define $H_{-1}$ norms based on the 
standard inner product $\langle\cdot,\cdot\rangle$:
\begin{equation}\label{equiv_norm}
\|\phi\|_{-1,\lambda} =\langle \phi, (\lambda-L)^{-1} \phi\rangle.
\end{equation}From \cite{S} we have that 
\begin{equation}
\alpha^{-1}\| \phi\|_{-1, \beta^{-1}\lambda}^{TASEP}  \le
\| \phi \|_{-1,\lambda} \le \alpha
\| \phi\|_{-1,\beta\lambda}^{TASEP} 
\end{equation}
  From the translation invariance of the generators 
\begin{eqnarray*}
\lll \phi \lll^2_{-1,\lambda}&=&\langle \sum_x \tau_x \phi , (\lambda-L)^{-1} \phi\rangle\\&=&\lim_{n\to\infty} \frac1{2n} \langle \sum_{x=-n}^n \tau_x \phi , (\lambda-L)^{-1} \sum_{x=-n}^n \tau_x\phi\rangle\\
&=&\lim_{n\to\infty} \frac1{2n}  \| \sum_{x=-n}^n \tau_x \phi \|^2_{-1,\lambda}.
\end{eqnarray*}
The proposition follows.\end{proof}

\begin{prop}\label{prop2}
Let $w$ be the current corresponding to a general AEP as in (\ref{w}) and $w^{TASEP}$ be the 
current for TASEP.
Then there exists $C<\infty$ such that
for  $0<\lambda<1$,
\begin{equation} 
\lll w-b w^{TASEP}\lll_{-1,\lambda} \le C.
\end{equation}
\end{prop}

{\it Remarks.}  1. In the theorem one can use either $L$
or $L^{TASEP}$ to define 
$\lll\cdot\lll_{-1,\lambda}$ since the results are equivalent.

2. This is similar to, but not the same as, 
results in \cite{SX}, because of the special 
norm $\lll\cdot\lll_{-1,\lambda}$.

\begin{proof} Since 
\[w-b w^{TASEP}=\sum_x x p(x) \left( \hat \eta_{\{0,1\}}-\hat \eta_{\{0,x\}}\right)\]
it is enough to show that
\begin{equation}\label{eq_sym}
\lll\hat \eta_{\{0,1\}}-\hat \eta_{\{0,k+1\}}\lll_{-1,\lambda}\le C
\end{equation}
for each $k>0$ where $C$ is a constant depending on $p(\cdot)$ and $k$ and $\lll\cdot\lll_{-1,\lambda}$ is defined using the generator  $L^{TASEP}$.
Call $V=\hat \eta_{\{0,1\}}-\hat \eta_{\{0,k+1\}}$.  
Dropping the third term in  the variational formula (\ref{varform}) we
 have
\begin{equation}
\lll V\lll^2_{-1,\lambda} \le \lla V, (\lambda-S)^{-1} V\rra.
\end{equation}
We now show that the right hand side is bounded independent of $\lambda$.  The computation is done using the fact that $S$
maps $H_2$ to itself. In particular, if $f,g\in H_2$ with
$
f=\sum_{x<y} f(x,y) \hat \eta_{\{x,y\}}$ and $ g=\sum_{x<y} g(x,y) \hat \eta_{\{x,y\}}$ then
\begin{equation}\label{reduction}
\lla f, g\rra=\sum_z \sum_{x<y} f(x+z,y+z) g(x,y)=\sum_{x=0}^\infty \overline{f}(x) \overline{g}(x).
\end{equation}
where 
\begin{equation}
\overline{f}(x)=\sum_y f(y,y+x+1),
\end{equation}
and  $Sf = \sum_{x<y} \hat S f(x,y) \hat \eta_{\{x,y\}}$ with
\begin{eqnarray}
\hat Sf(x,y)&=&\frac12\big(f(x,y+1)+f(x-1,y)-2 f(x,y)\\ \nonumber
&&\phantom{\frac12(}+\mathbf{1}_{\{y-x>1\}}\left(f(x,y-1)+f(x+1,y)-2 f(x,y)\right)\big).
\end{eqnarray}
Moreover
\begin{equation}\label{ssep_dual}
\overline{\hat Sf}(x)=(\mathcal{S}\bar f)(x)=\overline{f}(x+1)-\overline{f}(x)+\mathbf{1}_{\{x>0\}} (\overline{f}(x-1)- \overline{f}(x)).
\end{equation}
Our $V=\sum_{x<y} V(x,y) \hat \eta_{\{x,y\}}$ where $V(x,y)=\mathbf{1}_{\{x=0,y=1\}}-\mathbf{1}_{\{x=0,y=k+1\}}$. We can write
\begin{equation}
(\lambda-S)^{-1} V=\sum_{x<y} h(x,y) \hat \eta_{x,y}.
\end{equation}
for some $h$.
Then 
\begin{eqnarray}
 \lla V, (\lambda-S)^{-1} V\rra
 & = & \sum_x h(x,x+1)-\sum_x h(x,x+k+1)\nonumber
 \\
 & = & \overline{h}(0)-\overline{h}(k)
 \\\nonumber
 &=& (\lambda-\mathcal{S})^{-1} \overline{V}(0)-(\lambda-\mathcal{S})^{-1} \overline{V}(k)
\end{eqnarray}
where $\overline{V}(x)=\mathbf{1}_{\{x=0\}}-\mathbf{1}_{\{x=k\}}$.

An explicit computation shows that
\begin{equation}
q(x):=\frac{\gamma^x}{\lambda+1-\gamma}=\left((\lambda-\mathcal{S})^{-1} \mathbf{1}_{\{x=0\}}\right)(x)
\end{equation}
where $\gamma=\gamma(\lambda)$ is the solution of the equation
\begin{equation}
\lambda+2=\gamma^{-1}+\gamma
\end{equation}
with $0<\gamma<1$.
This is easy to check: if $x>0$ then
\[
\big((\lambda-S)q \big)(x)=\left(\lambda-(\gamma-1)-(\gamma^{-1}-1)\right)\frac{\gamma^x}{\lambda+1-\gamma}=0
\]
and 
\[
\big((\lambda-S)q\big)(0)=(\lambda-(\gamma-1))\frac{1}{\lambda+1-\gamma}=1.
\]
A similar calculation shows that if $k>0$ then one can find constants $c_1, c_2$ (depending on $k$ and $\lambda$) such that
\begin{eqnarray*}
\left((\lambda-\mathcal{S})^{-1} \mathbf{1}_{\{x=k\}}\right)(x)=
\left\{
\begin{array}{ccc}
\frac12(c_1 \gamma^{k-x}+c_2 \gamma^{x-k})&\textup{if}& 0\le x <k\\[5pt]
\frac12(c_1+c_2) \, \gamma^{x-k}&\textup{if}& k\le x
\end{array}
\right.
\end{eqnarray*}
and that  there is a $C<\infty$ such that
\begin{equation}
|c_i-\lambda^{-1/2}|\le C, \quad i=1,2. 
\end{equation}
So
\begin{equation}
(\lambda-\mathcal{S})^{-1} \overline{V}(0)-(\lambda-\mathcal{S})^{-1} \overline{V}(k)=\frac{1-\gamma^k}{\lambda+1-\gamma}+\frac12(c_1 (1-\gamma^k)+c_2 (1-\gamma^{-k})).
\end{equation}
Since $\gamma\simeq 1-\lambda^{1/2}$ as $\lambda\to 0$, it is not hard to check that the right hand side is bounded
for $0<\lambda<1$.
\end{proof}

\section{Monotonicity of $tD(t)$}

Let $X(t)$ be the position of a second class particle at time $t$ started at the origin and $
\widetilde X(t)=X(t)-(1-2\rho) b t$.

\begin{lemma}\label{prop:sym}
For any AEP,
\begin{equation}\label{iden1}
(1-\rho)E[\widetilde X(t)|\eta_{y}(0)=0]+ \rho E[\widetilde X(t)|\eta_{y}(0)=1]=0 
\end{equation}
and
\begin{equation}\label{iden2}
E[\widetilde X(t)|\eta_{y}(0)=1]=E[\widetilde X(t)|\eta_{-y}(0)=1]
\end{equation}
\end{lemma}
\begin{proof}
(\ref{iden1}) is straightforward from $E[X(t)]=(1-2\rho) t b$. To prove (\ref{iden2}) we first write the difference as 
\begin{eqnarray}
 &&\hskip-10mm \sum_x x(P(X(t)=x|\eta_{y}(0)=1)-P(X(t)=x|\eta_{-y}(0)=1))
\\
&&=\sum_x x P(X(t)=x|\eta_{y}(0)=1)-\sum_x (x+y) P(X(t)=x|\eta_{-y}(0)=1)+y\nonumber
\end{eqnarray}
We can write 
\begin{eqnarray}
&& \hskip-10mmP(X(t)=x|\eta_{y}(0)=1) \\ && = E[\eta_x(t)|\eta_0(0)=1,\eta_y(0)=1] -E[\eta_x(t)|\eta_0(0)=0,\eta_y(0)=1]\nonumber
\end{eqnarray}
and  by the translation invariance
\begin{eqnarray}
&& \hskip-10mmP(X(t)=x|\eta_{-y}(0)=1)\label{321}\\&&=E[\eta_x(t)|\eta_0(0)=1,\eta_{-y}(0)=1]-E[\eta_x(t)|\eta_0(0)=0,\eta_{-y}(0)=1]
\nonumber\\ \nonumber
&&=E[\eta_{x+y}(t)|\eta_0(0)=1,\eta_{y}(0)=1]-E[\eta_{x+y}(t)|\eta_0(0)=1,\eta_{y}(0)=0].
\end{eqnarray}
Substituting these into the previous equation we get
\begin{eqnarray}
&&\hskip-10mmE[\widetilde X(t)|\eta_{y}(0)=1]-E[\widetilde X(t)|\eta_{-y}(0)=1]\\\nonumber &&=\sum_x x \{E[\eta_x(t)|\eta_0(0)=1,\eta_{y}(0)=0] - E[\eta_x(t)|\eta_0(0)=0,\eta_{y}(0)=1]\}
+y\nonumber\\\nonumber
&&=\chi^{-1}\sum_x x E[\eta_x(t) (\eta_0(0)-\eta_y(0)]+y\\\nonumber &&=0
\end{eqnarray}
by (\ref{MV}).  
\end{proof}

\begin{lemma} For any AEP,
\begin{equation}
\partial_t (t D(t))=\sum_z z^2p(z)-2\rho\sum_{z>0}  z (p(z)-p(-z)) 
E[\widetilde X(t)|\eta_{z}(0)=1].
\end{equation}\label{lemma2}
\end{lemma}

\begin{proof}
We compute 
\begin{equation}\label{eq:sum324}
\partial_t \sum_x x^2 S(x,t)=\sum_{x,z} x^2   p(z) \la\eta_{x-z}(t) (1-\eta_{x}(t))-\eta_{x}(t) (1-\eta_{x+z}(t)),\eta_0(0)\ra.
\end{equation}
Summing by parts,  using the translation invariance, reversing space and time, we can rewrite (\ref{eq:sum324}) as
\begin{equation}\label{eq:325}
\sum_{x,z} (-2xz+z^2) p(-z) \la\eta_{0} (0)(1-\eta_{z}(0)),\eta_{x}(t)-\rho\ra\end{equation}
Again, by explicit computation $\la\eta_{0} (1-\eta_{z})(0),\eta_{x}(t)-\rho\ra$ is given by
\begin{equation}\label{326}
\chi \left(
\rho(1-\rho)\left(e_{11}-e_{01}\right)
+(1-\rho)^2\left(e_{10}-e_{00}\right)
-\rho \left(e_{11}-e_{10}\right)
\right)\end{equation}
where
\begin{equation}
e_{ij} = 
E[\eta_x(t)|\eta_0(0)=i,\eta_z(0)=j].\end{equation}
(\ref{326}) can be rewritten in terms of the second class particle (see (\ref{321})) as
\begin{equation}
\chi \left(
 (1-\rho) P(X(t)=x)- \rho P(X(t)=x-z|\eta_{-z}(0)=1)
\right).
\end{equation}
Substituting this into (\ref{eq:325}) and using (\ref{drift}), (\ref{MV}):
\begin{eqnarray}
\partial_t \sum_x x^2 S(x,t)&=&\chi (1-\rho) \sum_z (-2 E[X(t)] z+z^2) p(-z)\nonumber\\
&&-\chi \rho \sum_z (-2 E[X(t)|\eta_{-z}(0)=1] z-z^2) p(-z)\nonumber\\
&=&\chi \sum_z z^2 p(z)+2 b^2 t (1-2\rho)^2\label{eq:sok} \\&&\nonumber\hskip10mm-2 \chi \rho \sum_z E[\widetilde X(t)|\eta_{z}(0)=1]  z p(z)
\end{eqnarray}
Using (\ref{iden2}) and the definition of $D(t)$ completes the proof.
\end{proof}

\begin{prop}\label{prop:mon}
Suppose that  $p(z)=0$ for $|z|\neq 1$ (nearest neighbor). Then $t D(t)$ is 
 non-decreasing in $t$.
\end{prop}

\begin{proof}  We can assume
 $p(1)\ge p(-1)$.  In this case we will show
\begin{eqnarray}
E[\widetilde X(t)|\eta_{1}(0)=1]\le (1-\rho).
\end{eqnarray} 
By the previous lemma,
\begin{equation}
\partial_t (t D(t))\ge (1-2\rho(1-\rho)) p(1)
+(1+2\rho(1-\rho)) p(-1)\ge 0.
\end{equation}

 Consider a configuration where at site $1$ we have a second class particle, at site $0$ we have a third class particle and at all the other sites the distribution of particles is independent Bernoulli with probability $\rho$. The ordinary particles don't see the second or third class particles (i.e.~they see them as empty sites) and the second class particle doesn't see the third class particle. Let the process evolve according to the AEP dynamics, and denote the position of the third and second class particle with $A(t)$ and $B(t)$, respectively. It is not hard to see that the law of $A(t)$ is the same as the law of $X(t)$ conditioned on the event $\{\eta_1(0)=1\}$ and we have to prove
\begin{equation}
E[A(t)]\le (1-\rho) +(1-2 \rho) b t.
\end{equation}
Also, the law of $B(t)$ is the same as the law of $X(t)+1$ conditioned on the event $\{\eta_{-1}(0)=0\}$. By Lemma \ref{prop:sym} we have 
\[
E[\rho A(t)+(1-\rho) B(t)]=(1-\rho) +(1-2\rho) b t.
\]
thus it is enough to prove that
\begin{equation}\label{eq:aux}
E[A(t)] \le E[B(t)].
\end{equation}
Define the variable $Z(t)$ the following way: $Z(t)=0$ if $A(t)<B(t)$ and $Z(t)=1$ otherwise. Consider a possible joint trajectory $(x_1(t),x_2(t))$ for  $$\big(\min(A(t), B(t)),\max(A(t),B(t))\big).$$  Conditioned on $\{(x_1(s),x_2(s)), 0\le s \le t\}$, $Z(t)$  is a continuous time Markov process on $\{0,1\}$ 
with rate $p(-1){\mathbf 1}_{ \{x_2(t)-x_1(t)=1\}}$ for the transition $0\to 1$ and $p(1){\mathbf 1}_{ \{x_2(t)-x_1(t)=1\}}$ for the transition $1\to 0$.
This uses the fact that our process is nearest neighbor, and thus $Z(t)$ can change only if the second and third class particles switch places. We can now calculate $P(Z(t)=0\big| (x_1(s),x_2(s)), 0\le s \le t)$ explicitly. Let $T(t)=\left|\{s: x_2(s)-x_1(s)=1, 0\le s\le t \}\right|$ be the time spent by the two particles up to time $t$ with distance 1 between them, then (using $P(Z(0)=0)=1$)
\begin{equation}
P(Z(t)=0\big| (x_1(s),x_2(s)), 0\le s \le t)=\frac{p(-1) e^{-T(t) (p(-1)+p(1))} +p(1)}{p(-1)+p(1)}
\end{equation} 
Since  $p(1)\ge p(-1)$, this is always at least $1/2$. This means
\[
E\left[A(t)\big|\{(x_1(s),x_2(s)),0\le s\le t\}\right]\le E\left[B(t)\big|\{(x_1(s),x_2(s)),0\le s\le t\}\right]
\]
from which (\ref{eq:aux}) and the proposition follows. 
\end{proof}

%


\section{Summation by parts}\label{section_ps}

\newcommand{\sgn}{\operatorname{sgn}}

In this section we will prove identities (\ref{d-var}) and (\ref{abs_mom}).   They hold for general finite range exclusions, but we only need them in case of the TASEP so we will only give the proofs in that special case.  Note $b=1$ here.
The identities are a simple consequence of (\ref{stov}) and summation by parts,
once one knows the precise behaviour of $v(t,x)$ as $|x|\to \infty$.
They are not new; see, for example \cite{FF} for a proof of (\ref{abs_mom}).  But we could not find a reference for (\ref{d-var}),
so we include the proof here.

For $x\in \Z, t\ge 0$ denote by $N_t(x)$ the number 
of jumps from site $x$ to site $x+1$ up to time $t$ .

\begin{lemma} \label{lemmadvar}
\begin{eqnarray*}
v(x,t)&=&4 \chi |x|+ 4 Cov(N_t(0), N_t(x))-4 \sgn(x)\,  Cov\left(N_t(0), \sum\nolimits_{y=-|x|+1}^{|x|} \eta_y (t) \right).
\end{eqnarray*}

\end{lemma}


\begin{proof}
We will assume $x\ge 0$; the case $
x<0$ is analogous.
Recalling the definition (\ref{def_M}) of $M_t(x)$ and $N_t(x)$ we have
\begin{equation}\label{current}
N_t(0)-N_t(x)=\frac12 (M_t(x)-M_0(x)).
\end{equation}
It is easy to compute $Var(M_t(x))=4 \chi |x|$, and   
by the definition of $v(x,t)$ we have
\begin{eqnarray}
 v(x,t) \nonumber
&=&Var(2 N_t(0)-M_t(x))\\
&=& 4 Var(N_t(0))+4 \chi x-4 Cov(N_t(0),M_t(x))\label{term_}
\end{eqnarray}
Using the  identity (\ref{current}) and the translation invariance  we get
\begin{eqnarray}
 Cov(N_t(0),M_t(x))&=& E[N_t(0) M_t(x)]-E[N_t(0)] E[M_t(x)]
\\  \nonumber
&=& E[(N_t(x)+\frac12 (M_t(x)-M_0(x)) )M_t(x)]-E[N_t(0)] E[M_t(x)]\\
\nonumber
&=&E[(\frac12 (M_t(x)-M_0(x)))^2]+E[N_t(x) M_t(x)]-E[N_t(x)] E[M_t(x)]\\
\nonumber
&=&E[(N_t(x)-N_t(0))^2]+E[N_t(x) M_t(x)]-E[N_t(x)] E[M_t(x)]\\ \nonumber
&=&2 Var(N_t(0))-2Cov(N_t(x), N_t(0))+Cov(N_t(x), M_t(x)).
\end{eqnarray}
We will substitute this into (\ref{term_}) to get
\begin{eqnarray}
 v(x,t)&=& 4 \chi x+ 4Cov(N_t(x), N_t(0))\\\nonumber
 &&-2 Cov(N_t(0),M_t(x))-2  Cov(N_t(x),M_t(x))
\end{eqnarray}
By translation invariance, and because of the sign convention in the definition (\ref{def_M}) of $M_t(x)$, 
 \begin{equation}Cov(N_t(x), M_t(x))=-Cov(N_t(0), M_t(-x)),
 \end{equation}
and the lemma follows.
\end{proof}

\begin{lemma} \label{lem_cor}
For each $t>0$, there exist $C_1<\infty$ and $C_2>0$ such that
\begin{equation}
Cov\left(N_t(0), \eta_x (t) \right)\le C_1\exp\{-C_2 |x|\}, \quad Cov(N_t(0), N_t(x))\le C_1\exp\{-C_2 |x|\}.
\end{equation}
\end{lemma}

\begin{proof}
The lemma is standard, but we could not find an exact reference, so for completeness, we give a sketch of the proof.
Consider two copies $(\eta(t), \tilde\eta(t))$ of TASEP, coupled as in the preamble to (\ref{gra}), starting with initial data $(\eta_y, \tilde\eta_y=\eta_y \mathbf{1}_{\{y\in [-x/3,x/3]\cup[2x/3,4x/3]\}}), ~y\in \Bbb Z$ where $\eta$ is distributed according to $\pi_\rho$. Discrepancies perform nearest neighbor random walks, and the rate of jumping left or right is always at most 1. 
Let
\begin{equation}
A=\left\{\eta_0(s)=\tilde\eta_0(s)\quad \textup{and}\quad \eta_x(s)=\tilde\eta_x(s) \quad\textup{for~all}\quad s\in[0,t]\right\}. 
\end{equation}
$A^C$ is contained in the event that an initial discrepancy reaches $0$ or $x$ during the time interval $[0,t]$.  Because of the preservation of order, there are just $4$ candidates and hence $P(A^C) \le 4P(\textup{Poisson}(t)>x/3)$, which is exponentially small in $x$.
On $A$, $N_t(0)=\tilde N_t(0)$ and $N_t(x)=\tilde N_t(x)$, where $\tilde N_t(\cdot)$ are the currents in $\tilde\eta(t)$. Both $N_t(x)$ and $\tilde N_t(x)$ are stochastically dominated by Poisson($t$) random variables and hence, for any fixed $t$, their moments are bounded. Breaking up the respective expectations onto
$A$ and $A^C$ and applying Schwartz's inequality we see that both $|Cov(N_t(0), \eta_x (t) )-Cov(\tilde N_t(0), \tilde\eta_x (t) )|$ and $|Cov(N_t(0), N_t(x))-Cov(\tilde N_t(0), \tilde N_t(x))|$ are exponentially small in $x$. 

Hence it suffices to prove the lemma for the second process. Consider a
third process $\bar \eta$ with the same initial conditions as $\tilde \eta$, but disallowing jumps between $[2x/3]$ and $[2x/3]-1$. Using the
same argument as above,  it is enough to prove the lemma for $\bar\eta$. Now $\bar N_t(0)$ and $\bar \eta_x(t)$ (and $\bar N_t(0)$ and $\bar N_t(x)$) are independent, so the covariances vanish.
\end{proof}

Once one has Lemma \ref{lem_cor}, it follows from Lemma \ref{lemmadvar} that for each fixed $t\ge 0$, 
\begin{equation}\label{obv}
|v(x,t) -4\chi|x|| \le C_3\exp\{-C_4|x|\}
\end{equation}
 for some $C_3<\infty$ and $C_4>0$.  
Now (\ref{d-var}) and (\ref{abs_mom}) follow 
by taking partial summations, applying (\ref{stov})  
summing by parts, and noting that the boundary terms are exponentially
small from (\ref{obv}). 
\medskip

\noindent{\textbf{Acknowledgment.} The authors would like to thank the referee for pointing out an error in an earlier version of the manuscript.


\end{document}